\documentclass{article}

\usepackage{amscd}
\usepackage{amsfonts}
\usepackage{hyperref}
\usepackage{url}
\usepackage{xspace}

\usepackage{fancyvrb}
\usepackage{moreverb}
\usepackage{amsfonts}
\usepackage{helvet}
\usepackage{fancyheadings}
\usepackage{array, rotating} 
\usepackage{graphicx, amsmath}
\usepackage{texdraw}

 \DeclareMathOperator{\supp}{\rm supp} \DeclareMathOperator{\ddiv}{div}
\DeclareMathOperator{\eval}{eval} \DeclareMathOperator{\Div}{Div} 
\DeclareMathOperator{\Aut}{Aut}  \DeclareMathOperator{\Perm}{Perm}

\newtheorem{theorem}{Theorem}
\newtheorem{corollary}[theorem]{Corollary}

\newtheorem{lemma}[theorem]{Lemma}
\newtheorem{proposition}[theorem]{Proposition}

\newtheorem{example}[theorem]{Example}

\newtheorem{remark}{Remark}

\def\ppp{{\mathbb{P}}}

\def\zzz{\mathbb{Z}}

\def\oF{\overline{F}}
\def\pf{{\bf proof}:\ }
\def\qed{$\Box$}

\begin{document}

\title{Automorphism groups of generalized Reed-Solomon codes}

\author{David Joyner, Amy Ksir, Will Traves
\thanks{
Mathematics Dept., U.S. Naval Academy, Annapolis, MD 21402.
E-mails: wdj@usna.edu, ksir@usna.edu, traves@usna.edu}
}


\maketitle

\begin{abstract}
We look at AG codes associated to $\ppp^1$, re-examining the problem of determining their automorphism groups
(originally investigated by D\"ur in 1987 using combinatorial techniques) using recent methods from algebraic
geometry. We classify those finite groups that can arise as the automorphism group of an AG code and give an
explicit description of how these groups appear. We give examples of generalized Reed-Solomon codes with
large automorphism groups $G$, such as $G=PSL(2,q)$, and explicitly describe their $G$-module structure.
\end{abstract}



\section{Introduction}

Reed-Solomon codes are popular in applications because fast encoding and decoding algorithms are known for
them. For example, they are used in compact discs (more details can be found in \S 5.6 in Huffman and
Pless \cite{HP}).

In this paper we study which groups can arise as automorphism groups of a related collection of codes, the
algebraic geometry (AG) codes on $\ppp^1$. These codes are monomially equivalent to generalized Reed-Solomon
(GRS) codes. Their automorphism groups were first studied by D\"ur \cite{D} in 1987 using combinatorial
techniques. Huffman \cite{Hu} gives an excellent exposition of D\"ur's original work. In this paper, using
recent methods from algebraic geometry (due to Brandt and Stichenoth \cite{Sti2}, Valentini and
Madan \cite{VM}, Kontogeorgis \cite{K}), we present a method for computing GRS codes with ``large'' permutation
automorphism groups. In contrast to D\"ur's results, we indicate exactly how these automorphism groups can be
obtained.

The paper is organized as follows.  In section \ref{sec:codebackground} we review some background on AG codes
and GRS codes.  In section \ref{sec:action} we review some known results on automorphisms of AG codes, and
then prove our main result, characterizing the automorphism groups of AG codes. In section \ref{sec:examples}
we use these results to give examples of codes with large automorphism groups.  In section
\ref{sec:repstructure}, we discuss the structure of these group representations as $G$-modules, in some cases
determining it explicitly.

\section{AG codes and GRS codes}\label{sec:codebackground}

We recall some well-known background on AG codes and GRS codes.

Let $X$ be a smooth projective curve over a field $F$ and let $\oF$ denote a separable algebraic closure of
$F$.   We will generally take $F$ to be finite of order $q$.  Let $F(X)$ denote the function field of $X$
(the field of rational functions on $X$). Recall that a \textbf{divisor} on $X$ is a formal sum, with integer
coefficients, of places of $F(X)$. We will denote the group of divisors on $X$ by $\Div(X)$. The rational
points of $X$ are the places of degree $1$, and the set of rational points is denoted $X(F)$.

AG codes associated to a divisor $D$ are constructed from the Riemann-Roch space
\[
L(D)=L_X(D)= \{f\in F(X)^\times \ |\ \ddiv(f)+D\geq 0\}\cup \{0\},
\]
where $\ddiv(f)$ denotes the (principal) divisor of the function $f\in F(X)$. The Riemann-Roch space is a
finite dimensional vector space over $F$, whose dimension is given by the Riemann-Roch theorem.

Let $P_1,...,P_n\in X(F)$ be distinct points and $E=P_1+...+P_n\in \Div(X)$. Assume these divisors have
disjoint support, ${\rm supp}(D)\cap {\rm supp}(E)=\emptyset$. Let $C(D,E)$ denote the AG code

\begin{equation}
\label{eqn:AGcode} C(D,E) =\{(f(P_1),...,f(P_n))\ |\ f\in L(D)\}.
\end{equation}
This is the image of $L(D)$ under the evaluation map

\begin{equation}
\label{eqn:eval}
\begin{array}{c}
\eval_E:L(D)\rightarrow F^n,\\
f \longmapsto (f(P_1),...,f(P_n)).
\end{array}
\end{equation}
The following is well-known (a proof can be found in Joyner and Ksir \cite{JK3}).

\begin{lemma}
If $\deg(D)>\deg(E)$ then $\eval_E$ is injective.
\end{lemma}

In this paper, we restrict to the case where $X$ is the projective line $\ppp^1$ over $F$. In this case, if
$\deg D \geq 0$ then $\dim L(D) = \deg D + 1$, and otherwise $\dim L(D) = 0$.  Thus we will be interested in
the case where $\deg D \geq 0$.

In the special case when $D$ is a positive integer multiple of the point at infinity, then this construction
gives a Reed-Solomon code. More generally,
\[
C=\{(\alpha_1 f(P_1),...,\alpha_n f(P_n))\ |\ f\in L(\ell \cdot \infty)\},
\]
is called a {\bf generalized Reed-Solomon code} (or GRS code), where $\alpha_1,...,\alpha_n$ is a fixed set
of non-zero elements in $F$ (called ``multipliers'').

In fact, for a more general $D$, this construction gives a code which is monomially equivalent to a GRS code,
and which furthermore is MDS (that is, $n+1=k+d$, where $n$ is the length, $k$ is the dimension, and $d$ is
the minimum distance of the code). We say that two codes $C,C'$ of length $n$ are \textbf{monomially
equivalent} if there is an element of the group of monomial matrices $Mon_n(F)$ -- those matrices with
precisely one non-zero entry in each row and column -- (acting in the natural way on $F^n$) sending $C$ to
$C'$ (as $F$-vector spaces). Here, the set

\begin{lemma}
\label{lemma:grs} Let $X=\ppp^1/F$, $D$ be any divisor of positive degree on $X$, and let $E = P_1 + \cdots +
P_n$, where $P_1, \ldots P_n$ are points in $X(F)$ and $n > \deg D$.  Let $C(D,E)$ be the AG code constructed
as above.  Then $C(D,E)$ is an MDS code which is monomially equivalent to a GRS code (with all scalars
$\alpha_i=1$).
\end{lemma}

\pf  This is well-known (see for example Stichtenoth \cite{Sti}, \S II.2), but we give the details for
convenience. $C(D,E)$ has length $n$ and dimension $k = \deg(D)+1$.  By Theorem 13.4.3 of Huffman and
Pless \cite{HP}, its minimum distance $d$ satisfies
\[
n-\deg(D)\leq d,
\]
and the Singleton bound says that
\[
d \leq n+1-k = n - \deg(D).
\]
Therefore, $d = n+1 - k$, and this shows that $C(D,E)$ is MDS.

The  monomial equivalence follows from the fact that on $\ppp^1$, all divisors of a given positive degree are
(rationally) equivalent, so $D$ is rationally equivalent to $\deg(D) \cdot \infty$. Thus there is a rational
function $h$ on $X$ such that
\[
D = \deg(D) \cdot \infty + \ddiv(h).
\]
Then for any $f \in L(D)$, $fh$ is in $L(\deg(D) \cdot \infty)$.  Thus there is a map
\begin{eqnarray*}
M:  C(D,E) &\to&  C(\deg(D) \cdot \infty, E) \\
 (f(P_1), \ldots, f(P_n)) & \mapsto & (fh(P_1), \ldots, fh(P_n))
\end{eqnarray*}
which is linear and whose matrix is diagonal with diagonal entries $h(P_1), \ldots, h(P_n)$. In particular,
$M$ is a monomial matrix, so $C(D,E)$ and the GRS code $C(\deg(D) \cdot \infty, E)$ are monomially
equivalent. \qed

\begin{remark}
\label{rem:self-dual} The \textbf{spectrum} of a code of length $n$ is the list $[A_0,A_1,...,A_n]$, where
$A_i$ denotes the number of codewords of weight $i$. The \textbf{dual code} of a linear code $C\subset F^n$
is the dual of $C$ as a vector space with respect to the Hamming inner product on $F^n$, denoted $C^\perp$.
We say $C$ is \textbf{formally self-dual} if the spectrum of $C^\perp$ is the same as that of $C$. The
spectrum of any MDS code is known (see \S 7.4 in Huffman and Pless \cite{HP}), and as a consequence of this we
have the following
\[
A_j = \left(
\begin{array}{c}
n\\
j\end{array} \right) (q-1)\sum_{i=0}^{j-d}(-1)^i \left(
\begin{array}{c}
j-1\\
i\end{array} \right) q^{j-d-i},\ \ \ \ d\leq j\leq n,
\]
where $q$ is the order of the finite field $F$. The following is an easy consequence of this and the fact
that the dual code of an MDS code is MDS: if $C$ is a GRS code with parameters $[n,k,d]$ satisfying $n=2k$
then $C$ is formally self-dual.  We will see later some examples of formally self-dual codes with large
automorphism groups.
\end{remark}

\section{Automorphisms}\label{sec:action}
The action of a finite group $G\subset \Aut(X)$ on $F(X)$ is defined by restriction to $G$ of the map
\[
\begin{array}{cccc}
\rho:&\Aut(X)&\longrightarrow &\Aut(F(X)),\\
 & g &\longmapsto & (f\longmapsto f^g)
\end{array}
\]
where $f^g(x)=(\rho(g)(f))(x)=f(g^{-1}(x))$.

Note that $Y=X/G$ is also smooth and the quotient map

\begin{equation}
\label{eqn:psi} \psi:X\rightarrow Y
\end{equation}
yields an identification $F(Y)=F(X)^G:=\{f\in F(X)\ |\ f^g=f,\ \forall\, g\in G\}$.

Of course, $G$ also acts on the group $\Div(X)$ of divisors of $X$. If $g\in \Aut(X)$ and $d_P \in \zzz$, for
places $P$ of $F(X)$, then $g(\sum_P d_P P)=\sum_P d_Pg(P)$. It is easy to show that
$\ddiv(f^g)=g(\ddiv(f))$. Because of this, if $\ddiv(f)+D\geq 0$ then $\ddiv(f^g)+g(D)\geq 0$, for all $g\in
\Aut(X)$. In particular, if the action of $G$ on $X$ leaves $D\in \Div(X)$ stable then $G$ also acts on
$L(D)$. We denote this action by
\[
\rho:G\rightarrow \Aut(L(D)).
\]

Now suppose that $E = P_1 + \cdots + P_n$ is also stabilized by $G$.  In other words, $G$ acts on the set
$\supp(E) = \{P_1,  \ldots , P_n \}$ by permutation.  Then the group $G$ acts on $C(D,E)$ by $g\in G$ sending
$c=(f(P_1),...,f(P_n))\in C$ to $c'=(f(g^{-1}(P_1)),...,f(g^{-1}(P_n)))$, where $f\in L(D)$.

\begin{remark}
Observe that this map sending $c\longmapsto c'$, denoted $\phi(g)$, is well-defined. This is clearly true if
$\eval_E$ is injective. In case $\eval_E$ is not injective, suppose $c$ is also represented by $f'\in L(D)$,
so $c=(f'(P_1),...,f'(P_n))\in C$. Since $G$ acts on the set $\supp(E)$ by permutation, for each $P_i$,
$g^{-1}(P_i)=P_j$ for some $j$. Then $f(g^{-1}(P_i)) = f(P_j) = f'(P_j) = f'(g^{-1}(P_i))$, so
$(f(g^{-1}(P_1)),...,f(g^{-1}(P_n)))= (f'(g^{-1}(P_1)),...,f'(g^{-1}(P_n)))$. Therefore, $\phi(g)$ is
well-defined.
\end{remark}

The \textbf{permutation automorphism group} of the code $C$, denoted $\Perm(C)$, is the subgroup of the
symmetric group $S_n$ (acting in the natural way on $F^n$) which preserves the set of codewords. More
generally, we say two codes $C$ and $C'$ of length $n$ are \textbf{permutation equivalent} if there is an
element of $S_n$ sending $C$ to $C'$ (as $F$-vector spaces). The \textbf{automorphism group} of the code $C$,
denoted $\Aut(C)$, is the subgroup of the group of monomial matrices $Mon_n(F)$ (acting in the natural way on
$F^n$) which preserves the set of codewords. Thus the permutation automorphism group of $C$ is a subgroup of
the full automorphism group.

The map $\phi$ induces a homomorphism of $G$ into the automorphism group of the code. The image of the map
\begin{equation}
\label{eqn:phi}
\begin{array}{c}
\phi:G\rightarrow \Aut(C)\\
g\longmapsto \phi(g)
\end{array}
\end{equation}
is contained in $\Perm(C)$.

Define $\Aut_{D,E}(X)$ to be the subgroup of $\Aut(X)$ which preserves the divisors $D$ and $E$.

When does a group of permutation automorphisms of the code $C$ induce a group of automorphisms of the curve
$X$? Permutation automorphisms of the code $C(D,E)$ induce curve automorphisms whenever $D$ is very ample and
the degree of $E$ is large enough. Under these conditions, the groups $\Aut_{D,E}(X)$ and $\Perm C$ are
isomorphic.

\begin{theorem} (Joyner and Ksir \cite{JK2})
\label{thrm:lift} {\it Let $X$ be an algebraic curve, $D$ be a very ample divisor on $X$, and $P_1, \ldots,
P_n$ be a set of points on $X$ disjoint from the support of $D$.  Let $E = P_1 + \ldots + P_n$ be the
associated divisor, and $C = C(D,E)$ the associated AG code. Let $G$ be the group of permutation
automorphisms of $C$.  Then there is an integer $r \geq 1$ such that if $n > r \cdot \deg(D)$, then $G$ can
be lifted to a group of automorphisms of the curve $X$ itself. This lifting defines a group homomorphism
$\psi: \Perm C \to \Aut(X)$.  Furthermore, the lifted automorphisms will preserve $D$ and $E$, so the image
of $\psi$ will be contained in $\Aut_{D,E}(X)$. }
\end{theorem}

\begin{remark}
An explicit upper bound on $r$ can be determined (see Joyner-Ksir \cite{JK2}).  In the case where $X =
\ppp^1$, $r=2$. In addition, any divisor of positive degree on $\ppp^1$ is very ample.  Therefore, as long as
$\deg D >0$ and $n > 2 \deg (D)$, the groups $\Perm(C)$ and $\Aut_{D,E}(X)$ will be isomorphic.
\end{remark}

Now we would like to describe all possible finite groups of automorphisms of $\ppp^1$. Valentini and Madan
 \cite{VM} give a very explicit list of possible automorphisms of the associated function field $F(x)$ and
their ramifications.

\begin{proposition}
\label{prop:VM} (Valentini and Madan  \cite{VM}) Let $F$ be finite field of order $q=p^k$.  Let $G$ be a
nontrivial finite group of automorphisms of $F(x)$ fixing $F$ elementwise and let $E=F(x)^G$ be the fixed
field of $G$.  Let $r$ be the number of ramified places of $E$ in the extension $F(x)/E$ and $e_1,\ldots,e_r$
the corresponding ramification indices. Then $G$ is one of the following groups, with $F(x)/E$ having one of
the associated ramification behaviors:
\begin{enumerate}
\item
Cyclic group of order relatively prime to $p$ with $r = 2$, $e_1 = e_2 = | G_0|$.
\item
Dihedral group $D_m$ of order $2 m$ with $p = 2$, $( p, m) = 1$, $r = 2$, $e_1 = 2$, $e_2 = m$, or $p\not=
2$, $( p, m) = 1$, $r = 3$, $e_1 = e_2 = 2$, $e_3 = m$.
\item
Alternating group $A_4$ with $p\not= 2,3$, $r = 3$, $e_1 = 2$, $e_2 = e_3 = 3$.
\item
Symmetric group $S_4$ with $p\not= 2,3$, $r = 3$, $e_1 = 2$, $e_2 = 3$, $e_3 = 4$.
\item
Alternating group $A_5$ with $p = 3$, $r = 2$, $e_1 = 6$, $e_2 = 5$, or $p\not= 2, 3, 5$, $r = 3$, $e_1 = 2$,
$e_2 = 3$, $e_3 = 5$.
\item
Elementary Abelian $p$-group with $r = 1$, $e =|G_0|$.
\item
Semidirect product of an elementary Abelian $p$-group of order $q$ with a cyclic group of order $m$ with $m
|(q-1)$, $r = 2$, $e_1 = |G_0|$, $e_2 = m$.
\item
$PSL(2,q)$, with $p\not= 2$, $q = p^m$, $r = 2$, $e_1 = \frac{q(q-1)}{2}$, $e_2 = \frac{(q+1)}{2}$.
\item
$PGL(2,q)$, with $q = p^m$, $r = 2$, $e_1 = q(q-1)$, $e_2 = q+1$.
\end{enumerate}
\end{proposition}

The following result of Brandt can be found in \S 4 of Kontogeorgis and Antoniadis \cite{KA}. It provides a
more detailed explanation of the group action on $\ppp^1$ than the previous Proposition, giving the orbits
explicitly in each case.

Notation: In the result below, let $i=\sqrt{-1}$.

\begin{proposition}
\label{prop:B} (Brandt \cite{Br}) If the characteristic $p$ of the algebraically closed field of constants $F$
is zero or $p > 5$ then the possible automorphism groups of the projective line are given by the following
list.
\begin{enumerate}
\item
Cyclic group of order $\delta$.
\item
$D_\delta = \langle \sigma,\tau\rangle$, $(\delta , p ) = 1$ where $\sigma (x ) = \xi x$ , $\tau (x ) = 1/x$,
$\xi$ is a primitive $\delta$-th root of one. The possible orbits of the $D_\delta$ action are $B_\infty =
\{0, \infty\}$, $B^- = \{{\rm roots\ of\ }x^\delta- 1\}$, $B_+ = \{{\rm roots\ of\ }x^\delta+ 1\}$, $B_a =
\{{\rm roots\ of\ }x^{2\delta}+x^\delta+ 1\}$, where $a \in F- \{\pm 2\}$.
\item
$A_4 = \langle \sigma,\mu \rangle$, $\sigma (x ) = - x$ , $\mu (x ) = i\frac{x+1}{x-1}$, $i^2 = -1$. The
possible orbits of the action are the following sets: $B_0 = \{0, \infty, \pm 1, \pm i \}$, $B_1 = \{{\rm
roots\ of\ }x^4 - 2i\sqrt{3}x^2 + 1\}$, $B_2 = \{{\rm roots\ of\ }x^4 - 2i\sqrt{3}x^2 + 1\}$, $B_a = \{{\rm
roots\ of\ }\prod_{i =1}^3 (x^4 + a_i x^2 + 1)\}$, where $a_1 \in F- \{\pm 2, \pm 2i\sqrt{3}\}$, $a_2 =
\frac{2a_1+12}{2-a_1}$, $a_3 = \frac{2a_1-12}{2+a_1}$.
\item
$S_4 = \langle \sigma,\mu \rangle$, $\sigma (x ) = i x$, $\mu (x ) = i\frac{x+1}{x-1}$, $i^2 = -1$. The
possible orbits of the action are the following sets: $B_0 = \{0, \infty, \pm 1, \pm i \}$, $B_1 = \{{\rm
roots\ of\ }x^8 + 14x^4 + 1\}$, $B_2 = \{{\rm roots\ of\ }(x^4 + 1)(x^8 - 34x^4 + 1)\}$, $B_a = \{{\rm roots\
of\ }(x^8 + 14x^4 + 1)^3 - a (x^5 - x)^4 \}$, $a \in F- \{ 108\}$.
\item
$A_5 = \langle \sigma,\rho \rangle$, $\sigma (x ) = \xi x$, $\mu (x ) = -\frac{x+b}{bx+1}$, where $\xi$ is a
primitive fifth root of one and $b = -i (\xi^4 + \xi )$, $i^2 = -1$. The possible orbits of the action are
the following sets: $B_\infty = \{0, \infty\} \cup \{{\rm roots\ of\ }f_0 (x ) := x^{10} + 11i x^5 + 1\}$,
$B_0 = \{{\rm roots\ of\ } f_1 (x ) := x^{20} - 228i x^{15} - 494x^{10} - 228i x^5 + 1\}$, $B_0^* = \{{\rm
roots\ of\ } x^{30} + 522ix^{25} + 10005x^{20}- 10005x^{10} - 522ix^5 - 1\}$, $B_a = \{{\rm roots\ of\ } f_1
( x )^3 - a f_0 ( x )^5 \}$, where $a \in F - \{-1728i \}$.
\item
Semidirect products of elementary Abelian groups with cyclic groups: $(\zzz/p\zzz \times ... \times
\zzz/p\zzz)\times \zzz/m\zzz$ of order $p^t m$, where $m |(p^t - 1)$. Suppose we have an embedding of a field
of order $p^t$ into $k$. Assume $GF(p^t)$ contains all the $m$-th roots of unity. The possible orbits of the
action are the following sets: $B_\infty = \{\infty\}$, $B_0 = \{{\rm roots\ of\ } f (x ) := x\prod_{j
=1}^{(p^t-1)/m} (x^m - b_j )\}$, where $b_j$ are selected so that all the elements of the additive group
$\zzz/p\zzz \times ... \times \zzz/p\zzz$ ($t$ times), when viewed as elements in $F$, are roots of $f (x )$,
$B_a = \{{\rm roots\ of\ }f(x)^m - a \}$, where $a \in F - B_0$.
\item
$PSL(2, p^t ) = \langle \sigma,\tau, \phi \rangle$, $\sigma (x ) = \xi^2 x$ , $\tau (x ) = -1/x$,
$\phi(x)=x+1$, where $\xi$ is a primitive $m=p^t-1$ root of one. The orbits of the action are $B_\infty =
\{\infty, {\rm roots\ of\ }x^m-x\}$. $B_0 = \{{\rm roots\ of\ } (x^m - x)^{m -1} + 1\}$, $B_a = \{{\rm roots\
of\ } ((x^m - x)^{m -1} + 1)^{(m+1)/2}-a(x^m - x)^{m(m -1)/2}\}$, where $a \in F^\times$.
\item
$PGL(2, p^t ) = \langle \sigma,\tau, \phi \rangle$, $\sigma (x ) = \xi x$ , $\tau (x ) = 1/x$, $\phi(x)=x+1$,
where $\xi$ is a primitive $m=p^t-1$ root of one. The orbits of the action are $B_\infty = \{\infty, {\rm
roots\ of\ }x^m-x\}$. $B_0 = \{{\rm roots\ of\ } (x^m - x)^{m -1} + 1\}$, $B_a = \{{\rm roots\ of\ } ((x^m -
x)^{m -1} + 1)^{m+1}-a(x^m - x)^{m(m -1)}\}$, where $a \in F^\times$.
\end{enumerate}
\end{proposition}

\pf  Brandt \cite{Br}, Stichtenoth \cite{Sti2}.

\qed

Let $Y = X/G$ be the curve associated to the field $E$ in Proposition \ref{prop:VM}, and let $\pi: X \to Y$
be the quotient map.

\begin{corollary}
\label{cor:main} Assume that (1) the finite field $F$ has characteristic $>5$, (2) $\pi$ is defined over $F$,
(3) for each $p_1\in X(F)$, all the points $p_0$ in the fiber $\pi^{-1}(p_1)$ are rational: $p_0\in X(F)$,
and (4) $F$ is so large that the orbits described in Proposition \ref{prop:B} are complete.  Then the above
Proposition \ref{prop:B} holds over $F$.
\end{corollary}

\pf  Under the hypotheses given, the inertia group is always equal to the decomposition group and the action
of the group $G$ of automorphisms commutes with the action of the absolute Galois group $\Gamma = {\rm
Gal}(\oF/F)$.

\qed

The following is our main result.

\begin{theorem}
{\it Assume $C$ is a GRS code constructed from a divisor $D$ with positive degree and defined over a
sufficiently large finite field $F$ (as described in Corollary \ref{cor:main}). Then the automorphism group
of $C$ must be one of the groups in Proposition \ref{prop:VM}. }
\end{theorem}

In fact, the action can be made explicit using Proposition \ref{prop:B}.

\begin{corollary}
Each GRS code over a sufficiently large finite field is monomially equivalent to a code whose automorphism
group is one of the groups in Proposition \ref{prop:VM}.
\end{corollary}

\pf  (of theorem) We assume the field is as in Corollary \ref{cor:main}. Use Theorem \ref{thrm:lift} and Lemma
\ref{lemma:grs}.

\qed

It would be interesting to know if this result can be refined in the case when $n=2k$, as that might give
rise to a class of easily constructable self-dual codes with large automorphism group.

\section{Examples}\label{sec:examples}

Pick two distinct orbits ${\cal O}_1$ and ${\cal O}_2$ of $G$ in $X(F)$. Assume that $D$ is the sum of the
points in the orbit ${\cal O}_1$ and let ${\cal O}_2=\{P_1,...,P_n\}\subset X(F)$. Define the associated code
of length $n$ by

\[
C=\{(f(P_1),...,f(P_n))\ |\ f\in L(D)\}\subset F^n.
\]
This code has a $G$-action, by $g\in G$ sending $(f(P_1),...,f(P_n))$ to
\newline
$(f(g^{-1}P_1),...,f(g^{-1}P_n))$, so is a $G$-module. Indeed, by construction, the action of $G$ is by
permuting the coordinates of $C$.

\begin{example}
\label{ex:3} {\small{ {\it Let $F$ be a finite field of characteristic $>5$ which contains (1) all $4^{th}$
and $5^{th}$ roots of unity, (2) all the roots of $x^{10} + 11i x^5 + 1$, (3) all the roots $B_0$ of $x^{20}
- 228i x^{15} - 494x^{10} - 228i x^5 + 1$, and (4) all the roots $B_0^*$ of $x^{30} + 522ix^{25} +
10005x^{20}- 10005x^{10} - 522ix^5 - 1$. Furthermore, let $B_\infty = \{0, \infty\} \cup \{\text{roots of }
x^{10} + 11i x^5 + 1$\}. Let $E=\sum_{P\in B_0} P$ and let $D=\sum_{P\in B_0^*\cup B_\infty} P$. Then
$\deg(E)=20$ and $\deg(D)=42$. Then $C=C(D,E)$ is a formally self-dual code with parameters $n=42$, $k=21$,
$d=22$, and automorphism group $A_5$. } }}
\end{example}

This follows from (5) of Proposition \ref{prop:B} and Remark \ref{rem:self-dual}.

\begin{example}
\label{ex:4} {\small{ {\it Let $F=GF(q)$ be a finite field of characteristic $p>5$ for which $q\equiv 1\pmod
8$ and for which $F$ contains (1) all the roots of $x^{q-1} - x$, and (2) all the roots $B_1$ of
$((x^{q-1}-x)^{q-2}+1)^{q/2} - (x^{q-1}-x)^{(q-1)(q-2)/2}$. If $B_\infty = \{\infty, \text{roots of } x^{q-1}
- x\}$, then let $D=\frac{(q-1)(q-2)}{4}\sum_{P\in B_\infty} P$, $E=\sum_{P\in B_1} P$, and $C=C(D,E)$. Then
$C$ is a formally self-dual code with parameters $n=\frac{q(q-1)(q-2)}{2}$, $k=n/2$, $d=n+1-k$, and
permutation automorphism group $G=PSL(2,q)$. } }}
\end{example}

This follows from (7) of Proposition \ref{prop:B}.

\section{Structure of the representations}
\label{sec:repstructure}

We study the possible representations of finite groups $G$ on the codes $C(D,E)$. As noted in Lemma
\ref{thrm:lift}, when $E$ is large enough, this is the same as the representation of $G$ on $L(D)$.
Therefore we study the possible representations of $G$ on $L(D)$.  For simplicity we will restrict to the
case where the support of $D$ is rational, i.e. $D=\sum_{i=1}^s a_i P_i$, where $P_1, \ldots, P_s$ are
rational points on $\ppp^1$.

We can give the representation explicitly by finding a basis for $L(D)$. For a divisor $D$ with rational
support on $X = \ppp^1$, it is easy to find a basis for $L(D)$, as follows. Let $\infty=[1:0]\in X$ denote
the point corresponding to the localization $\overline{F}[x]_{(1/x)}$, and $[p:1]$ denote the point
corresponding to the localization $\overline{F}[x]_{(x-p)}$, for $p\in \overline{F}$. For notational
simplicity, let
\[
m_P(x)= \left\{
\begin{array}{cc}
x, & P=[1:0]=\infty,\\
\frac{1}{(x-p)}, & P=[p:1].
\end{array}
\right.
\]
Then $m_P(x)$ is a rational function with a simple pole at the point $P$, and no other poles.

\begin{lemma}
\label{lemma:P1basis} Let $D=\sum_{i=1}^s a_i P_i$ be a divisor with rational support on $X=\ppp^1$, so
$a_i\in\zzz$ and $P_i \in X(F)$ for $0\leq i\leq s$.
\begin{itemize}
\item[(a)]
If $D$ is effective then
\[
\{ 1,m_{P_i}(x)^k\ |\ 1\leq k\leq a_i, 1\leq i\leq s\}
\]
is a basis for $L(D)$.

\item[(b)]
If $D$ is not effective but $\deg(D)\geq 0$ then $D$ can be written as $D=D_1+D_2$, where $D_1$ is effective
and $\deg(D_2)=0$.  Let $q(x)\in L(D_2)$ (which is a 1-dimensional vector space) be any non-zero element.
Let $D_1=\sum_{i=1}^s a_i P_i$. Then
\[
\{ q(x), m_{P_i}(x)^{k}q(x)\ |\ 1\leq k\leq a_i, 1 \leq i \leq s \}
\]
is a basis for $L(D)$.

\item[(c)]
If $\deg(D)<0$ then $L(D)=\{0\}$.
\end{itemize}

\end{lemma}


\pf   This is an easy application of the Riemann-Roch theorem. Note that the first part appears as Lemma 2.4
in Lorenzini \cite{L}.

By the Riemann-Roch theorem, $L(D)$ has dimension $\deg D +1$ if $\deg (D) \geq 0$ and otherwise $L(D) =
\{0\}$, proving part (c) and the existence of $q(x)$ in part (b). For part (a), since $m_{P_i}(x)^k$ has a
pole of order $k$ at $P_i$ and no other poles, it will be in $L(D)$ if and only if $k \leq a_i$. Similarly,
for part (b), $m_{P_i}(x)^k$ will be in $L(D_1)$ if and only if $k \leq a_i$; therefore $m_{P_i}(x)^k q(x)$
will be in $L(D_1+D_2)=L(D)$ under the same conditions.  In each of parts (a) and (b), the set of functions
given is linearly independent, so by a dimension count must form a basis for $L(D)$.  \qed

Now let $G$ be a finite group acting on $X=\ppp^1$ and let $D$ be a divisor with rational support, stabilized
by $G$.  Let $S=\supp(D)$ and let
\[
S=S_1\cup S_2\cup ... \cup S_m
\]
be the decomposition of $S$ into primitive $G$-sets.  Then we can write $D$ as
\[
D=\sum_{k=1}^m a_k S_k = \sum_{k=1}^m a_k \sum_{i=1}^s P_{ik},
\]
where for each $k$, $P_{1k} \ldots P_{sk}$ are the points in the orbit $S_k$. Then $G$ will act by a
permutation on the points $P_{1k} \ldots P_{sk}$ in each orbit, and therefore on the corresponding functions
$m_{P_{sk}}(x)$.

\begin{theorem}
\label{thrm:P1} Let $X$, $F$, $G\subset \Aut(X)=PGL(2,\overline{F})$, and $D$ be as above. Let
$\rho:G\rightarrow \Aut(L(D))$ denote the associated representation.

\begin{itemize}
\item[(a)]
If $D$ is effective then
\[
\rho\cong \mathbf{1} \oplus_{k=1}^m a_k \rho_k,
\]
where $\mathbf{1}$ denotes the trivial representation, and $\rho_k$ is the permutation representation on the
subspace
\[
V_k= \mbox{ span } \{ m_P(x)\ |\ P\in S_k\}.
\]

\item[(b)]
If $\deg(D)>0$ but $D$ is not effective then $L(D)$ is a sub-$G$-module of $L(D^+)$, where $D^+$ is a
$G$-invariant effective divisor satisfying $D^+\geq D$.
\end{itemize}

\end{theorem}


The groups and orbits which can arise are described in Proposition \ref{prop:VM} above.

\pf  (a) By part (a) of Lemma \ref{lemma:P1basis}, $\{1,m_{P_{ik}}(x)^\ell\ |\ 1\leq \ell\leq a_i, 1\leq i\leq
s\, 1 \leq k \leq m\}$ form a basis for $L(D)$.  $G$ will act trivially on the constants.  For each $\ell$,
$G$ will act by permutations as described on each set $\{ m_{P_{ik}}(x)^{\ell}\ |\ P_{ik}\in S_k\}.$

(b) Since $D$ is not effective, we may write $D=D^+-D^-$, where $D^+$ and $D^-$ are non-zero effective
divisors. The action of $G$ must preserve $D^+$ and $D^-$. Since $L(D)$ is a $G$-submodule of $L(D^+)$, the
claim follows.

\qed

\vskip .1in

{\it Acknowledgements}: We thank Cary Huffman for very useful suggestions 
on an earlier version and for the
references to D\"ur \cite{D} and Huffman \cite{Hu}.
We also thank John Little for valuable suggestions that
improved the exposition.

\end{document}